\documentclass[a4paper,11pt]{amsart}
\usepackage{graphicx}
\usepackage{mathptmx}
\usepackage{amsmath}
\usepackage{amssymb}
\usepackage{enumitem}
\usepackage{xcolor}
\usepackage{xparse}
\NewDocumentCommand{\eulerian}{omm}
 {%
  \genfrac<>{0pt}{}{#2}{#3}%
  \IfValueT{#1}{_{\!#1}}%
 }

 \newmuskip\pFqmuskip

\newcommand*\pFq[6][8]{%
  \begingroup 
  \pFqmuskip=#1mu\relax
  \mathchardef\normalcomma=\mathcode`,
  \mathcode`\,=\string"8000
  \begingroup\lccode`\~=`\,
  \lowercase{\endgroup\let~}\pFqcomma
  {}_{#2}F_{#3}{\left(\genfrac..{0pt}{}{#4}{#5}\bigg|#6\right)}%
  \endgroup
}
\newcommand{\pFqcomma}{{\normalcomma}\mskip\pFqmuskip}

\newtheorem{theorem}{Theorem}

\newtheorem{corollary}[theorem]{Corollary}

\begin{document}

\title[Some formulas for fully degenerate Bernoulli numbers and polynomials]{Some formulas for fully degenerate Bernoulli numbers and polynomials}

\author{Taekyun  Kim}
\address{Department of Mathematics, Kwangwoon University, Seoul 139-701, Republic of Korea}
\email{tkkim@kw.ac.kr}

\author{DAE SAN KIM}
\address{Department of Mathematics, Sogang University, Seoul 121-742, Republic of Korea}
\email{dskim@sogang.ac.kr}



\maketitle

\begin{abstract}
The aim of this paper is to study the fully degenerate Bernoulli polynomials and numbers, which are a degenerate version of  Bernoulli polynomials and numbers and arise naturally from the Volkenborn integral of the degenerate exponential functions on $\mathbb{Z}_p$.
We find some explicit expressions for the fully degenerate Bernoulli polynomials and numbers in terms of the degenerate Stirling numbers of the second kind, the degenerate $r$-Stirling numbers of the second kind and of the degenerate Stirling polynomials. We also consider the degenerate poly-Bernoulli polynomials and derive explicit representations for them in terms of the same degenerate Stirling numbers and polynomials.
\end{abstract}

\section{Introduction}

In recent years, studies on degenerate versions of special polynomials and numbers, which started with the paper of Carlitz in [4], have regained interests of some mathematicians. They have been done with various tools, including combinatorial methods, generating functions, umbral calculus, $p$-adic analysis, probability theory, differential equations, operator theory, analytic number theory and special functions. \par
The aim of this paper is to study the fully degenerate Bernoulli polynomials and numbers (see \eqref{3}), which are a degenerate version of Bernoulli polynomials and numbers. They are different from  Carlitz degenerate Bernoulli polynomials and numbers (see \eqref{2}) and arise naturally from the Volkenborn integral of the degenerate exponential functions on $\mathbb{Z}_p$ (see \eqref{3-1}).
We find some explicit expressions for the fully degenerate Bernoulli polynomials and numbers in terms of the degenerate Stirling numbers of the second kind, the degenerate $r$-Stirling numbers of the second kind and of the degenerate Stirling polynomials (see \eqref{4}, \eqref{5}, \eqref{34-1}). We also consider the degenerate poly-Bernoulli polynomials and derive explicit representations for them in terms of the same degenerate Stirling numbers and polynomials. The degenerate two variable Fubini polynomials (see \eqref{9}) are degenerate versions of the two variable Fubini polynomials. We note here that the generating function of the fully degenerate Bernoulli polynomials is the integral over the unit interval of the generating function of the degenerate two variable Fubini polynomials. This fact is used repeatedly throughout this paper. For the rest of this section, we recall some necessary facts that will be used throughout this paper.

\vspace{0.1in}

For any $\lambda\in\mathbb{R}$, the degenerate exponential functions are defined by
\begin{equation}
e_{\lambda}^{x}(t)=\sum_{k=0}^{\infty}\frac{(x)_{k,\lambda}}{k!}t^{k},\quad (\mathrm{see}\ [11,15,22]).\label{1}	
\end{equation}
where $(x)_{0,\lambda}=1,\ (x)_{n,\lambda}=x(x-\lambda)\cdots(x-(n-1)\lambda),\ (n\ge 1)$. Note that $\lim_{\lambda\rightarrow 0}e_{\lambda}^{x}(t)=e^{xt}$.\par
For $x=1$, we use notation $e_{\lambda}(t)=e_{\lambda}^{1}(t)$. \par
In [4], Carlitz introduced the degenerate Bernoulli polynomials given by
\begin{equation}
\frac{t}{e_{\lambda}(t)-1}e_{\lambda}^{x}(t)=\sum_{n=0}^{\infty}\beta_{n}(x|\lambda)\frac{t^{n}}{n!}.\label{2}
\end{equation}
When $x=0$, $\beta_{n}(\lambda)=\beta_{n}(0|\lambda)$ are called the degenerate Bernoulli numbers. \par
In 2016, a variant of the Carlitz degenerate Bernoulli polynomials, called the fully degenerate Bernoulli polynomials, is defined by
\begin{equation}
\frac{\log(1+\lambda t)}{\lambda(e_{\lambda}(t)-1)}e_{\lambda}^{x}(t)=\sum_{n=0}^{\infty}\beta_{n,\lambda}(x)\frac{t^{n}}{n!},\quad (\mathrm{see}\ [12,19]).\label{3}
\end{equation}
When $x=0$, $\beta_{n,\lambda}=\beta_{n,\lambda}(0)$ are called the fully degenerate numbers. \par
Note that $\displaystyle\lim_{\lambda\rightarrow 0}\beta_{n,\lambda}(x)=B_{n}(x),\ (n\ge 0),\displaystyle$ where $B_{n}(x)$ are the ordinary Bernoulli polynomials given by
\begin{displaymath}
	\frac{t}{e^{t}-1}e^{xt}=\sum_{n=0}^{\infty}B_{n}(x)\frac {t^{n}}{n!},\quad (\mathrm{see}\ [1-23]).
\end{displaymath}

Before proceeding further, we would like to explain that $\beta_{n,\lambda}(x)$ arise naturally as the Volkenborn integral of the degenerate exponential in \eqref{1} on $\mathbb{Z}_p$. Let $p$ be a fixed prime number. Let $\mathbb{Z}_{p},\mathbb{Q}_{p}$ and $\mathbb{C}_{p}$ denote the ring of $p$-adic integers, the field of $p$-adic  rational numbers and the completion of the algebraic closure of $\mathbb{Q}_{p}$, respectively. Here the $p$-adic norm $|\cdot|_{p}$ is normalized as $|p|_{p}=\frac{1}{p}$. \par
For any $\mathbb{C}_p$-valued uniformly differential function $f$ on $\mathbb{Z}_{p}$, the Volkenborn integral of $f$ on $\mathbb{Z}_p$ is defined by
\begin{align*}
\int_{\mathbb{Z}_{p}}f(x)d\mu(x) = \lim_{n\rightarrow\infty}\frac{1}{p^{n}}\sum_{x=0}^{p^{n}-1}f(x),
\end{align*}
which satisfies the fundamental property
\begin{align*}
\int_{\mathbb{Z}_{p}}f(x+1)d\mu(x)- \int_{\mathbb{Z}_{p}}f(x)d\mu(x)=f^{\prime}(0).
\end{align*}
By using this property, for $\lambda,t \in\mathbb{C}_{p}$ with $|\lambda t|_{p}<p^{-\frac{1}{p-1}}$, we have
\begin{align}
\int_{\mathbb{Z}_{p}}e_{\lambda}^{x+y}(t)d\mu(y)=\frac{\log(1+\lambda t)}{\lambda(e_{\lambda}(t)-1)}e_{\lambda}^{x}(t).\label{3-1}
\end{align}
Then we get from \eqref{1} and \eqref{3-1} the following Witt type identity for $\beta_{n,\lambda}(x)$:
\begin{align*}
\beta_{n,\lambda}(x)=\int_{\mathbb{Z}_p}(x+y)_{n,\lambda}d\mu(y),
\end{align*}
This shows that the fully degenerate Bernoulli polynomials and numbers indeed arise naturally as $p$-adic integrals on $\mathbb{Z}_p$. \par
It is known that the degenerate Stirling numbers of the second kind are defined by
\begin{equation}
(x)_{n,\lambda}=\sum_{k=0}^{n}S_{2,\lambda}(n,k)(x)_{k},\quad (n\ge 0),\quad (\mathrm{see}\ [12,13,19]),\label{4}
\end{equation}
where $(x)_{0}=1,\ (x)_{n}=x(x-1)\cdots(x-n+1),\ (n\ge 1)$. \par
For $r\in\mathbb{Z}$ with $r\ge 0$, the degenerate $r$-Stirling numbers of the second kind are defined by
\begin{equation}
(x+r)_{n,\lambda}=\sum_{k=0}^{n}S_{2,\lambda}^{(r)}(n+r,k+r)(x)_{k},\ \ (n\ge 0),\quad (\mathrm{see}\ [18,20]).\label{5}
\end{equation}
From \eqref{4} and \eqref{5}, we note that
\begin{equation}
\frac{1}{k!}\Big(e_{\lambda}(t)-1\Big)^{k}=\sum_{n=k}^{\infty}S_{2,\lambda}(n,k)\frac{t^{n}}{n!},\quad (\mathrm{see}\ [11]),\label{6}	
\end{equation}
and
\begin{equation}
e_{\lambda}^{r}(t)\frac{1}{k!}\Big(e_{\lambda}(t)-1\Big)^{k}=\sum_{n=k}^{\infty}S_{2,\lambda}^{(r)}(n+r,k+r)\frac{t^{n}}{n!},\quad (\mathrm{see}\ [11,20]).\label{7}	
\end{equation}
Note that $\displaystyle\lim_{\lambda\rightarrow 0}S_{2,\lambda}(n,k)=S_{2}(n,k)\displaystyle$ and $\displaystyle\lim_{\lambda\rightarrow 0}S_{2,\lambda}^{(r)}(n,k)=S_{2}^{(r)}(n,k)\displaystyle$, where $S_{2}(n,k)$ and $S_{2}^{(r)}(n+r,k+r)$ are respectively the Stirling numbers of the second kind and the $r$-Stirling numbers of the second kind.
It is well known that two variable Fubini polynomials are defined by
\begin{equation}
\frac{1}{1-x(e^{t}-1)}e^{yt}=\sum_{n=0}^{\infty}F_{n}(x|y)\frac{t^{n}}{n!},\quad (\mathrm{see}\ [8,9]).\label{8}
\end{equation}
When $y=0$, $F_{n}(x|0)=F_{n}(x)$ are called the Fubini polynomials,. \par
From \eqref{8}, we note that
\begin{displaymath}
F_{n}(x|r)=\sum_{k=0}^{n}x^{k}k!S_{2}(n+r,k+r),\quad(\mathrm{see}\ [8,9]),	
\end{displaymath}
where $n,r\in\mathbb{Z}$ with $n,r\ge 0$.\par
In [10], the degenerate two variable Fubini polynomials are defined by
\begin{equation}
\frac{1}{1-x(e_{\lambda}(t)-1)}e_{\lambda}^{y}(t)=\sum_{n=0}^{\infty}F_{n,\lambda}(x|y)\frac{t^{n}}{n!}.\label{9}	
\end{equation}
Note that
\begin{equation}
F_{n,\lambda}(x|y)=\sum_{m=0}^{n}\sum_{k=0}^{m}\binom{n}{m}S_{2,\lambda}(m,k)k!x^{k}(y)_{n-m,\lambda},\quad (\mathrm{see}\ [10,16]),\label{10}
\end{equation}
where $n$ is a nonnegative integer. \par
When $y=0$, $F_{n,\lambda}(x)=F_{n,\lambda}(x|0)$ are called the degenerate Fubini polynomials (see [16]). \par

\section{Some formulas for fully degenerate Bernoulli numbers and polynomials}
From \eqref{6}, we have
\begin{align}
\sum_{n=k}^{\infty}S_{2,\lambda}(n,k)\frac{t^{n}}{n!}&=\frac{1}{k!}\Big(e_{\lambda}(t)-1\Big)^{k}=\frac{1}{k!}\sum_{j=0}^{k}\binom{k}{j}
(-1)^{k-j}e_{\lambda}^{j}(t) \label{11} \\
&=\sum_{n=0}^{\infty}\bigg(\frac{1}{k!}\sum_{j=0}^{k}\binom{k}{j}(-1)^{k-j}(j)_{n,\lambda}\bigg)\frac{t^{n}}{n!}.\nonumber
\end{align}
Thus, we get
\begin{equation}
S_{2,\lambda}(n,k)=\frac{1}{k!}\sum_{j=0}^{k}\binom{k}{j}(-1)^{k-j}(j)_{n,\lambda},\quad (n,k\ge 0).\label{12}
\end{equation}
From \eqref{7}, we note that
\begin{equation}
\sum_{n=k}^{\infty}S_{2,\lambda}^{(r)}(n+r,k+r)\frac{t^{n}}{n!}=e_{\lambda}^{r}(t)\frac{1}{k!}\Big(e_{\lambda}(t)-1\Big)^{k}=\sum_{n=0}^{\infty}\bigg(\frac{1}{k!}\sum_{j=0}^{k}\binom{k}{j}(-1)^{k-j}(r+j)_{n,\lambda}\bigg)\frac{t^{n}}{n!}.\label{13}	
\end{equation}
Comparing the coefficients on both sides of \eqref{13}, we have
\begin{equation}
S_{2,\lambda}^{(r)}(n+r,k+r)=\frac{1}{k!}\sum_{j=0}^{k}\binom{k}{j}(-1)^{k-j}(r+j)_{n,\lambda},\quad (n,k\ge 0).\label{14}	
\end{equation}
By \eqref{9}, we get
\begin{equation}
\frac{1}{1+y(e_{\lambda}(t)-1)}=\sum_{n=0}^{\infty}F_{n,\lambda}(-y)\frac{t^{n}}{n!}.\label{15}
\end{equation}
Thus, we have
\begin{equation}
F_{n,\lambda}(-y)=\sum_{k=0}^{n}S_{2,\lambda}(n,k)(-1)^{k}k!y^{k},\label{16}	
\end{equation}
and
\begin{equation}
\int_{0}^{x}F_{n,\lambda}(-y)dy=\sum_{k=0}^{n}S_{2,\lambda}(n,k)(-1)^{k}\frac{k!}{k+1}x^{k+1}.\label{17}
\end{equation}
Now, we observe that
\begin{align}
\frac{1}{e_{\lambda}(t)-1}\log\Big(1+x(e_{\lambda}(t)-1)\Big)&=\sum_{n=0}^{\infty}\int_{0}^{x}F_{n,\lambda}(-y)dy\frac{t^{n}}{n!}\label{18} \\
&=\sum_{n=0}^{\infty}\bigg(\sum_{k=0}^{n}S_{2,\lambda}(n,k)\frac{k!}{k+1}x^{k+1}(-1)^{k}\bigg)\frac{t^{n}}{n!}.\nonumber
\end{align}
Let us take $x=1$ in \eqref{18}. Then we have
\begin{equation}
\frac{1}{\lambda(e_{\lambda}(t)-1)}\log(1+\lambda t)=\sum_{n=0}^{\infty}\bigg(\sum_{k=0}^{n}S_{2,\lambda}(n,k)\frac{k!}{k+1}(-1)^{k}\bigg)\frac{t^{n}}{n!}.\label{19}	
\end{equation}
By \eqref{3} and \eqref{19}, we get
\begin{equation}
\beta_{n,\lambda}=\sum_{k=0}^{n}\frac{k!}{k+1}S_{2,\lambda}(n,k)(-1)^{k},\quad (n\ge 0).\label{20}
\end{equation}
From \eqref{12} and \eqref{20}, we have the following theorem.
\begin{theorem}
For $n\ge 0$, we have
\begin{displaymath}
\beta_{n,\lambda}=\sum_{k=0}^{n}\frac{(-1)^{k}}{k+1}k!S_{2,\lambda}(n,k)=\sum_{k=0}^{n}\frac{1}{k+1}\sum_{j=0}^{k}\binom{k}{j}(-1)^{j}(j)_{n,\lambda}.
\end{displaymath}	
\end{theorem}
Let $\triangle$ be the difference operator with $\triangle f(x)=f(x+1)-f(x)$. Then we have
\begin{equation}
\triangle^{k}f(x)=\sum_{j=0}^{k}\binom{k}{j}(-1)^{k-j}f(j+x),\quad (k\ge 0). \label{22}
\end{equation}
Let us take $f(x)=(x)_{n,\lambda},\ (n\ge 0)$. Then we have
\begin{align}
\triangle^{k}(0)_{n,\lambda}=\sum_{j=0}^{k}\binom{k}{j}(-1)^{k-j}(j)_{n,\lambda}=k!S_{2,\lambda}(n,k). \label{23}
\end{align}
Therefore, by Theorem 1 and \eqref{23}, we obtain the following corollary.
\begin{corollary}
For $n\ge 0$, we have
\begin{displaymath}
\beta_{n,\lambda}=\sum_{k=0}^{n}\frac{(-1)^{k}}{k+1}\triangle^{k}(0)_{n,\lambda}.
\end{displaymath}
\end{corollary}
Let $r$ be a nonnegative integer. Then, by \eqref{9}, we get
\begin{equation}
\frac{e_{\lambda}^{r}(t)}{e_{\lambda}(t)-1}\log\Big(1+x(e_{\lambda}(t)-1)\Big)=\sum_{n=0}^{\infty}\int_{0}^{x}F_{n,\lambda}(-y|r)dy\frac{t^{n}}{n!}.\label{24}
\end{equation}
From \eqref{7} and \eqref{9}, we note that
\begin{align}
\sum_{n=0}^{\infty}F_{n,\lambda}(-y|r)\frac{t^{n}}{n!}&=\frac{1}{1+y(e_{\lambda}(t)-1)}e_{\lambda}^{r}(t)\label{25} \\
&=\sum_{k=0}^{\infty}(-1)^{k}y^{k}k!\frac{1}{k!}\Big(e_{\lambda}(t)-1\Big)^{k}e_{\lambda}^{r}(t)\nonumber \\
&= \sum_{k=0}^{\infty}(-1)^{k}y^{k}k!\sum_{n=k}^{\infty}S_{2,\lambda}^{(r)}(n+r,k+r)\frac{t^{n}}{n!}\nonumber \\
&=\sum_{n=0}^{\infty}\bigg(\sum_{k=0}^{n}(-1)^{k}y^{k}k!S_{2,\lambda}^{(r)}(n+r,k+r)\bigg)\frac{t^{n}}{n!}.\nonumber	
\end{align}
By comparing the coefficients on both sides of \eqref{25}, we get
\begin{equation}
F_{n,\lambda}(-y|r)=\sum_{k=0}^{n}(-1)^{k}y^{k}k!S_{2,\lambda}^{(r)}(n+r,k+r).\label{26}
\end{equation}
From \eqref{24} and \eqref{26}, we can derive the following:
\begin{align}
&\frac{e_{\lambda}^{r}(t)}{e_{\lambda}(t)-1}\log\Big(1+x(e_{\lambda}(t)-1)\Big)=\sum_{n=0}^{\infty}\int_{0}^{x}F_{n,\lambda}(-y|r)dy\frac{t^{n}}{n}\label{27}\\
&=\sum_{n=0}^{\infty}\bigg(\sum_{k=0}^{n}(-1)^{k}\frac{x^{k+1}}{k+1}k!S_{2,\lambda}^{(r)}(n+r,k+r)\bigg)\frac{t^{n}}{n!}.\nonumber	
\end{align}
Let us take $x=1$ in \eqref{27}. Then, by \eqref{3}, \eqref{14} and \eqref{27}, we get
\begin{align}
\sum_{n=0}^{\infty}\beta_{n,\lambda}(r)\frac{t^{n}}{n!}&=\frac{\log(1+\lambda t)}{\lambda(e_{\lambda}(t)-1)}\cdot e_{\lambda}^{r}(t) \label{28} \\
&=\sum_{n=0}^{\infty}\bigg(\sum_{k=0}^{n}(-1)^{k}\frac{1}{k+1}k!S_{2,\lambda}^{(r)}(n+r,k+r)\bigg)\frac{t^{n}}{n!}\nonumber \\
&=\sum_{n=0}^{\infty}\bigg(\sum_{k=0}^{n}(-1)^{k}\frac{k!}{k+1}\frac{1}{k!}\sum_{j=0}^{k}(-1)^{k-j}\binom{k}{j}(r+j)_{n,\lambda}\bigg)\frac{t^{n}}{n!}.\nonumber
\end{align}
Now, the next theorem follows from \eqref{28}.
\begin{theorem}
For $n\ge 0$, we have
\begin{align*}
\beta_{n,\lambda}(r)=\sum_{k=0}^{n}(-1)^{k}\frac{1}{k+1}k!S_{2,\lambda}^{(r)}(n+r,k+r)=\sum_{k=0}^{n}\frac{1}{k+1}\sum_{j=0}^{k}\binom{k}{j}(-1)^{j}(r+j)_{n,\lambda}.
\end{align*}	
\end{theorem}
Let us take $f(x)=(x)_{n,\lambda},\ (n\ge 0),$ in \eqref{22}. Then we have
\begin{equation}
\triangle^{k}(r)_{n,\lambda}=\sum_{j=0}^{k}\binom{k}{j}(-1)^{k-j}(j+r)_{n,\lambda},\quad (n\ge 0).\label{30}
\end{equation}
Therefore, by \eqref{30}, we obtain the following corollary.
\begin{corollary}
	For $n\ge 0$, we have
	\begin{displaymath}
		\beta_{n,\lambda}(r)=\sum_{k=0}^{n}\frac{(-1)^{k}}{k+1}\triangle^{k}(r)_{n,\lambda}.
	\end{displaymath}
\end{corollary}
Now, we observe from \eqref{13} that
\begin{align}
\sum_{n=0}^{\infty}\frac{1}{k!}\sum_{j=0}^{k}\binom{k}{j}(-1)^{k-j}(j+r)_{n,\lambda}\frac{t^{n}}{n!}&=e_{\lambda}^{r}(t)\frac{1}{k!}\Big(e_{\lambda}(t)-1\Big)^{k}\label{31}\\
&=\sum_{m=0}^{\infty}(r)_{m,\lambda}	\frac{t^{m}}{m!}\sum_{l=k}^{\infty}S_{2,\lambda}(l,k)\frac{t^{l}}{l!}\nonumber \\
&=\sum_{n=k}^{\infty}\bigg(\sum_{l=k}^{n}\binom{n}{l}S_{2,\lambda}(l,k)(r)_{n-l,\lambda}\bigg)\frac{t^{n}}{n!}.\nonumber
\end{align}
By comparing the coefficients on both sides of \eqref{31}, we obtain the following theorem.
\begin{theorem}
For $k,n\ge 0$, we have
	\begin{align*}
\sum_{j=0}^{k}\binom{k}{j}(-1)^{j}(j+r)_{n,\lambda}&=(-1)^{k}k!\sum_{l=k}^{n}\binom{n}{l}S_{2,\lambda}(l,k)(r)_{n-l,\lambda}\\
&=(-1)^{k}\triangle^{k}(r)_{n,\lambda}.
	\end{align*}
\end{theorem}
From \eqref{1}, \eqref{6} and \eqref{7}, we note that
\begin{equation}
S_{2,\lambda}^{(r)}(n+r,k+r)=\sum_{l=k}^{n}\binom{n}{l}S_{2,\lambda}(l,k)(r)_{n-l,\lambda},\quad (n,k\ge 0),\label{32}
\end{equation}
and
\begin{equation}
(x+y)_{n,\lambda}=\sum_{l=0}^{n}\binom{n}{l}(x)_{l,\lambda}(y)_{n-l,\lambda},\quad (n\ge 0).\label{33}
\end{equation}
By \eqref{33}, we get
\begin{align}
\sum_{j=0}^{k}\binom{k}{j}(-1)^{j}(r+j)_{n,\lambda}&=\sum_{j=0}^{k}\binom{k}{j}(-1)^{j}\sum_{l=0}^{n}\binom{n}{l}(r)_{n-l,\lambda}(j)_{l,\lambda}\nonumber\\
&=\sum_{l=0}^{n}\binom{n}{l}(r)_{n-l,\lambda}\sum_{j=0}^{k}\binom{k}{j}(-1)^{j}(j)_{l,\lambda}\label{34} \\
&=\sum_{l=0}^{n}\binom{n}{l}(r)_{n-l,\lambda}S_{2,\lambda}(l,k)k!(-1)^{k}.\nonumber	
\end{align}
Thus, for $n,k\ge 0$ with $n\ge k$, we have
\begin{displaymath}
\sum_{j=0}^{k}\binom{k}{j}(-1)^{j}(r+j)_{n,\lambda}=\sum_{l=k}^{n}\binom{n}{l}(r)_{n-l,\lambda}(-1)^{k}k!S_{2,\lambda}(l,k).
\end{displaymath}
The degenerate Stirling polynomials are defined by
\begin{align}
\frac{1}{k!}\Big(e_{\lambda}(t)-1\Big)^{k}e_{\lambda}^{x}(t)=\sum_{n=k}^{\infty}S_{2,\lambda}(n,k|x)\frac{t^{n}}{n!},\quad (\mathrm{see}\ [14]) \label{34-1}
\end{align}
where $k$ is nonnegative integer. \par
Then, we easily get
\begin{align}
S_{2,\lambda}(n,k|x)&=\sum_{l=k}^{n}\binom{n}{l}S_{2,\lambda}(l,k)(x)_{n-l,\lambda},\quad (n\ge 0)\label{35}\\
&=\frac{1}{k!}\sum_{l=0}^{k}\binom{k}{l}(-1)^{k-l}(l+x)_{n,\lambda},\quad (\mathrm{see}\ [14]).\nonumber
\end{align}
From \eqref{25}, we note that
\begin{equation}
F_{n,\lambda}(-y|x)=\sum_{k=0}^{n}(-1)^{k}y^{k}k!S_{2,\lambda}(n,k|x).\label{36}	
\end{equation}
Thus, by \eqref{36}, we get
\begin{align}
\frac{e_{\lambda}^{x}(t)}{e_{\lambda}(t)-1}\log\Big(1+z(e_{\lambda}(t)-1)\Big)&=\sum_{n=0}^{\infty}\int_{0}^{z}F_{n,\lambda}(-y|x)dy\frac{t^{n}}{n!}\label{37}\\
&=\sum_{n=0}^{\infty}\bigg(\sum_{k=0}^{n}(-1)^{k}\frac{z^{k+1}}{k+1}k!S_{2,\lambda}(n,k|x)\bigg)\frac{t^{n}}{n!}.\nonumber	
\end{align}
Let us take $z=1$ in \eqref{37}. Then we have
\begin{align}
\sum_{n=0}^{\infty}\beta_{n,\lambda}(x)\frac{t^{n}}{n!}&=\frac{\log(1+\lambda t)}{\lambda(e_{\lambda}(t)-1)}e_{\lambda}^{x}(t)=\sum_{n=0}^{\infty}\sum_{k=0}^{n}\frac{(-1)^{k}}{k+1}k!S_{2,\lambda}(n,k|x)\frac{t^{n}}{n!}\label{38} \\
&=\sum_{n=0}^{\infty}\bigg(\sum_{k=0}^{n}(-1)^{k}\frac{k!}{k+1}\frac{1}{k!}\sum_{l=0}^{k}\binom{k}{l}(-1)^{k-l}(l+x)_{n,\lambda}\bigg)\frac{t^{n}}{n!}.\nonumber
\end{align}
Therefore, by \eqref{38}, we obtain the following theorem.
\begin{theorem}
For $n\ge 0$, we have
\begin{align*}
\beta_{n,\lambda}(x)=\sum_{k=0}^{n}\frac{(-1)^{k}}{k+1}k!S_{2,\lambda}(n,k|x)
=\sum_{k=0}^{n}\frac{1}{k+1}\sum_{l=0}^{k}\binom{k}{l}(-1)^{l}(l+x)_{n,\lambda}.
\end{align*}	
\end{theorem}
Note that
\begin{equation}
\triangle^{k}(x)_{n,\lambda}=\sum_{j=0}^{k}\binom{k}{j}(-1)^{k-j}(j+x)_{n,\lambda}.\label{39}
\end{equation}
Hence, by \eqref{39}, we get
\begin{displaymath}
	\beta_{n,\lambda}(x)=\sum_{k=0}^{n}\frac{(-1)^{k}}{k+1}\triangle^{k}(x)_{n,\lambda},\quad (n\ge 0).
\end{displaymath}
In [13,17,19], the degenerate polylogarithm function of index $k$ is defined by
\begin{equation}
\mathrm{Li}_{k,\lambda}(x)=\sum_{n=1}^{\infty}\frac{(-\lambda)^{n-1}(1)_{n,1/\lambda}}{(n-1)!n^{k}}x^{n},\quad (k\in\mathbb{Z}). \label{40}	
\end{equation}
Let $\log_{\lambda}t$ be the compositional inverse of $e_{\lambda}(t)$. Then we have
\begin{equation}
\log_{\lambda}(1+t)=\sum_{n=1}^{\infty}\lambda^{n-1}(1)_{n,1/\lambda}\frac{t^{n}}{n!},\quad (\mathrm{see}\ [11]). \label{41}	
\end{equation}
By \eqref{40} and \eqref{41}, we get
\begin{equation}
\mathrm{Li}_{1,\lambda}(t)=-\log_{\lambda}(1-t).\label{42}
\end{equation}
The degenerate poly-Bernoulli polynomials of index $k$ are defined by
\begin{equation}
\frac{\mathrm{Li}_{k,\lambda}(1-e_{\lambda}(-t))}{1-e_{\lambda}(-t)}e_{\lambda}^{-x}(-t)=\sum_{n=0}^{\infty}\beta_{n,\lambda}^{(k)}(x)\frac{t^{n}}{n!}.\label{43}
\end{equation}
Note that
\begin{align}
\frac{\mathrm{Li}_{1,\lambda}(1-e_{\lambda}(-t))}{1-e_{\lambda}(-t)}e_{\lambda}^{-x}(-t)&=\frac{t}{1-e_{\lambda}(-t)}e_{\lambda}^{-x}(-t)\label{44}\\
&=\frac{t}{e_{-\lambda}(t)-1}e_{-\lambda}^{(x+1)}(t)=\sum_{n=0}^{\infty}\beta_{n}(x+1|-\lambda)\frac{t^{n}}{n!},\nonumber	
\end{align}
where $\beta_{n}(x|\lambda)$ are the Carlitz degenerate Bernoulli polynomials in \eqref{2}. \par
From \eqref{43} and \eqref{44}, we have
\begin{equation}
\beta_{n,\lambda}^{(1)}(x)=\beta_{n}(x+1|-\lambda),\quad (n\ge 0).\label{45}	
\end{equation}
Now, we observe that
\begin{align}
\frac{\mathrm{Li}_{p,\lambda}(1-e_{\lambda}(-t))}{1-e_{\lambda}(-t)}e_{\lambda}^{-x}(-t)
&=\sum_{k=0}^{\infty}\frac{(-\lambda)^{k}(1)_{k+1,1/\lambda}}{(k+1)^{p}}\frac{1}{k!}\Big(e_{-\lambda}(t)-1\Big)^{k}e_{-\lambda}^{x-k}(t)\label{46}\\
&=\sum_{k=0}^{\infty}\frac{(-\lambda)^{k}(1)_{k+1,1/\lambda}}{(k+1)^{p}}\sum_{n=k}^{\infty}S_{2,-\lambda}(n,k|x-k)\frac{t^{n}}{n!}\nonumber \\
&=\sum_{n=0}^{\infty}\bigg(\sum_{k=0}^{n}\frac{(-\lambda)^{k}(1)_{k+1,1/\lambda}}{(k+1)^{p}}S_{2,-\lambda}(n,k|x-k)\bigg)\frac{t^{n}}{n!}.\nonumber
\end{align}
Therefore, by \eqref{43}, and \eqref{46}, we obtain the following theorem.
\begin{theorem}
For $p\in\mathbb{Z}$ and $n\ge 0$, we have
\begin{displaymath}
	\beta_{n,\lambda}^{(p)}(x)=\sum_{k=0}^{n}\frac{(-\lambda)^{k}(1)_{k+1,1/\lambda}}{(k+1)^{p}}S_{2,-\lambda}(n,k|x-k).
\end{displaymath}	
\end{theorem}
From \eqref{35}, we note that
\begin{align}
S_{2,-\lambda}(n,k|x-k)&=\frac{1}{k!}\sum_{l=0}^{k}\binom{k}{l}(-1)^{k-l}(l+x-k)_{n,-\lambda}\label{47}\\
&=\frac{1}{k!}\sum_{l=0}^{k}\binom{k}{l}(-1)^{l}(x-l)_{n,-\lambda}.\nonumber	
\end{align}
Therefore, by Theorem 7 and \eqref{47}, we obtain the following theorem.
\begin{theorem}
For $p\in\mathbb{Z}$ and $n\ge 0$, we have
\begin{displaymath}
\beta_{n,\lambda}^{(p)}(x)=\sum_{k=0}^{n}\frac{(-\lambda)^{k}(1)_{k+1,1/\lambda}}{(k+1)^{p}k!}\sum_{l=0}^{k}\binom{k}{l}(-1)^{l}(x-l)_{n,-\lambda}.
\end{displaymath}
\end{theorem}
Note that
\begin{align}
\sum_{l=0}^{k}\binom{k}{l}(-1)^{l}(x-l)_{n,-\lambda}&=\sum_{l=0}^{k}	\binom{k}{l}(-1)^{l}\sum_{j=0}^{n}\binom{n}{j}(-1)^{j}(x)_{n-j,-\lambda}(l)_{j,\lambda}\label{48}\\
&=\sum_{j=0}^{n}\binom{n}{j}(-1)^{j}(x)_{n-j,-\lambda}\sum_{l=0}^{k}\binom{k}{l}(-1)^{l}(l)_{j,\lambda}\nonumber \\
&=\sum_{j=0}^{n}\binom{n}{j}(-1)^{j}(x)_{n-j,-\lambda}k!(-1)^{k}S_{2,\lambda}(j,k).\nonumber
\end{align}
Therefore, by Theorem 8 and \eqref{48}, we obtain the following theorem.
\begin{theorem}
For $n\ge 0$ and $p\in\mathbb{Z}$, we have
\begin{displaymath}
	\beta_{n,\lambda}^{(p)}(x)=\sum_{k=0}^{n}\frac{(\lambda)^{k}(1)_{k+1,1/\lambda}}{(k+1)^{p}}\sum_{j=0}^{n}\binom{n}{j}(-1)^{j}S_{2,\lambda}(j,k)(x)_{n-j,-\lambda}.
\end{displaymath}	
\end{theorem}
From \eqref{14}, we note that
\begin{align}
\sum_{n=k}^{\infty}S_{2,\lambda}^{(r)}(n+r,k+r)\frac{t^{n}}{n!}&=e_{\lambda}^{r}(t)\frac{1}{k!}\Big(e_{\lambda}(t)-1\Big)^{k} \label{49} \\
&=\sum_{n=k}^{\infty}\bigg(\sum_{j=k}^{n}S_{2,\lambda}(j,k)(r)_{n-j,\lambda}\binom{n}{j}\bigg)\frac{t^{n}}{n!}.\nonumber	
\end{align}
For $r\in\mathbb{N}\cup\{0\}$, by Theorem 9 and \eqref{49}, we get
\begin{align}
\beta_{n,\lambda}^{(p)}(-r)&= \sum_{k=0}^{n}\frac{(\lambda)^{k}(1)_{k+1,1/\lambda}}{(k+1)^{p}}\sum_{j=0}^{n}\binom{n}{j}(-1)^{j}S_{2,\lambda}(j,k)(-r)_{n-j,-\lambda}\label{50} \\
&=(-1)^{n} \sum_{k=0}^{n}\frac{(\lambda)^{k}(1)_{k+1,1/\lambda}}{(k+1)^{p}}\sum_{j=0}^{n}\binom{n}{j}S_{2,\lambda}(j,k)(r)_{n-j,\lambda}\nonumber\\
&=(-1)^{n} \sum_{k=0}^{n}\frac{(\lambda)^{k}(1)_{k+1,1/\lambda}}{(k+1)^{p}}S_{2,\lambda}^{(r)}(n+r,k+r).\nonumber	
\end{align}
Therefore, by \eqref{50}, we obtain the following theorem.
\begin{theorem}
For $r,p\in\mathbb{Z}$ with $r\ge 0$, we have
\begin{align*}
\beta_{n,\lambda}^{(p)}(-r)&=(-1)^{n} \sum_{k=0}^{n}\frac{(\lambda)^{k}(1)_{k+1,1/\lambda}}{(k+1)^{p}}\sum_{j=0}^{n}\binom{n}{j}S_{2,\lambda}(j,k)(r)_{n-j,\lambda}\\
&=(-1)^{n} \sum_{k=0}^{n}\frac{(\lambda)^{k}(1)_{k+1,1/\lambda}}{(k+1)^{p}}S_{2,\lambda}^{(r)}(n+r,k+r).	
\end{align*}
\end{theorem}
\section{Conclusion}
In recent years, a lot of mathematicians have explored degenerate versions of many special numbers and polynomials, which were initiated by Carlitz. Their interests were not only in combinatorial and arithmetical properties but also in applications to differential equations, identities of symmetry and probability theory. These degenerate versions include the degenerate Stirling numbers of the first and second kinds, degenerate Bernoulli numbers of the second kind and degenerate Bell numbers and polynomials. \par
In this paper, we studied the fully degenerate Bernoulli polynomials and numbers that are a degenerate version of Bernoulli polynomials and numbers, and different from the Carlitz degenerate Bernoulli polynomials and numbers. We found several explicit expressions for the fully degenerate Bernoulli polynomials and numbers in terms of the degenerate Stirling numbers of the second, the degenerate $r$-Stirling numbers of the second kind and of the degenerate Stirling polynomials. We also investigated the degenerate poly-Bernoulli polynomials and derived explicit representations in terms of the same degenerate Stirling numbers and polynomials. \par
It is one of our future research projects to continue to study various degenerate versions of many special polynomials and numbers with applications to physics, science and engineering in mind.

\end{document}